\documentclass[11pt]{article}
\usepackage{graphicx}
\usepackage{enumerate}
\usepackage{amsmath}
\usepackage{amsfonts}
\usepackage{amssymb}
\usepackage{amsfonts}
\usepackage{amssymb,amsmath,amsthm,epsfig,graphicx,subfigure,hyperref}
\usepackage[margin=1.3in]{geometry}
\usepackage{caption2}
\usepackage{color}
\usepackage{enumerate}
\usepackage{booktabs}
\usepackage[square, sort, numbers]{natbib}
\usepackage{supertabular}

\RequirePackage{lineno}

\def\blackbox{\hfill {\vrule height3pt width4pt depth2pt}}
\def\Box{\hfill \framebox(5.25,5.25){}}

\newtheorem{thm}{Theorem}
\newtheorem{pro}{Proposition}
\newtheorem{lem}{Lemma}

\newtheorem{corr}{Corollary}
\newtheorem{defin}{Definition}
\newtheorem{ob}{Observation}

\newenvironment{observation}{\noindent \begin{ob}}{{\hfill$\Box$}\end{ob}}

\newenvironment{proposition}{\begin{pro} \nopagebreak}{{\hfill$\Box$} \end{pro}}
\newenvironment{thm-prf}{\begin{thm} \nopagebreak}{\end{thm}}
\newenvironment{pro-prf}{\begin{pro} \nopagebreak}{\end{pro}}

\newenvironment{lem-prf}{\begin{lem} \nopagebreak}{\end{lem}}

\newenvironment{corr-prf}{\begin{corr} \nopagebreak}{\end{corr}}

\newcommand {\be}{\begin{equation}}
\newcommand {\ee}{\end{equation}}
\begin{document}

\title{Integrated Facility Location and Production Scheduling in
Multi-Generation Energy Systems}
\author{Qiao-Chu~He and Tao~Hong \thanks{%
Q. He and T. Hong are with the Department of Systems Engineering and
Engineering Management, University of North Carolina, Charlotte, NC, 28223
USA. Email: qhe4@uncc.edu;hong@uncc.edu.} }
\maketitle

\begin{abstract}
In this paper, we investigate the energy system design problems with the 
\emph{multi-generation technologies}, i.e., simultaneous generation of
multiple types of energy. Our results illustrate the economic value of
multi-generation technologies to reduce spatio-temporal demand uncertainty
by risk pooling both within and across different facilities.

\end{abstract}

\newpage

\section{Introduction}

% and "HIS" in caps to complete the first word

Energy systems are multi-dimensional: Electricity, heat, oil, natural gas
and even biofuels, are distributed to meet demand with supply. An
interesting feature of the energy systems is that their dimensionalities are
intertwined: Multiple types of energy can be generated simultaneously;
Consequently, the optimal design for the value chain of one type of energy,
i.e., from generation to distribution, should not be independent from
another. This, however, poses great challenges to policy makers and system
planners.

In this paper, we investigate the energy systems design problems when
multi-generation technologies are available. Co-generation of heat and
electricity is a best-known example of such technologies, which increases
overall energy efficiency compared with separate production of heat and
power \cite{martens1998energetic}. By making more efficient use of fuel
inputs and renewable sources, multi-generation technologies also allow low
carbon emission \cite{meunier2002co}. A prominent trend in multi-generation
technologies is that, different energy solutions and devices are
increasingly installed at the users' premises to supply their local
multi-energy needs \cite{chicco2009distributed}.

Our research contributes to the state-of-the-art from a \emph{system
perspective}. In this paper, we address the following two research questions:

\begin{itemize}
\item How to design an energy production network which integrates
macro-level strategic decisions (facility location, multi-generation
technology investment, etc.) and micro-level operational decisions
(production planning, energy transportation etc.)?

\item What is the economic value of the multi-generation technologies under
demand uncertainties for energy?
\end{itemize}

To be specific, the first research question explores \emph{methodologies}
for the system design problems. The second research question is motivated by
the fact that both \emph{inter-temporal} and \emph{spatial} demand
uncertainties pose fundamental challenges to the design and optimization of
energy systems. For example, adequate cooling demand in the summertime and
thermal demand in the wintertime are needed to make their joint generation
economically feasible \cite{meunier2002co}. From this perspective, we not
only make methodological contributions to tackle this challenge, but also
offer insights on the value of the multi-generation technologies to
policy-makers.

The rest of this paper is organized as follows. Section \ref{s-lit} reviews
relevant literature. Section~\ref{s-model} introduces our model setup. In
Section \ref{s-analysis}, we carry out the analysis. In Section \ref{s-num},
we describe our computational method, and provide numerical examples.
Section~\ref{s-con} concludes this paper with a discussion of the future
research directions.

\section{Related Work}

\label{s-lit} Multi-generation technologies in power and energy systems
become increasingly attractive to the research community. A holistic
multi-energy system assessment (gas, heat and power) is conducted in \cite%
{883bfc133a82486082cff0a5ea0e86fe}. The flexibility design problem to cope
with the uncertainty in wind power generation is addressed in \cite%
{ma2011optimizing}. In addition, investment in system flexibility helps the
integration of wind power to the combined heat and power network \cite%
{chen2015increasing}. The existing literature often relies on variations of
the unit commitment model \cite%
{martinez2012unified,quelhas2007multiperiod,li2008interdependency,qiu2015low, martinez2013robust}%
, or dynamic control \cite{mancarella2013real}, to solve a dynamic
scheduling problem. However, this approach restricts its application to the
micro-level operational decisions due to its computational complexity. In
contrast to the aforementioned traditional approach, we model the production
scheduling component using an asymptotic \emph{proxy}, which implicitly
implements a dynamic allocation of production capacity, while abstracting
away from the scheduling details. This formulation is capable of being
integrated with both macro- and micro-level decisions in energy systems.

Facility location problems have been studied extensively in the operations
research community. The integrated approach of network design is proposed in 
\cite{Shen2003}, which is shown to be superior than solving location and
capacity planning problem separately. Our model is also related to the
literature on the transportation cost minimization \cite{Zipkin1980}, and
service scheduling problem \cite{Shanthikumar1992}. The integration of the
process flexibility into the network design has been studied by \cite%
{Mak2009a,Mak2009b}.

Our model is also related to network flexibility design literature, such as 
\cite{Deng2013} and the references therein. Nevertheless, we avoid any 
\textit{a priori} flexibility assumption by taking an engineering approach
(stochastic programming with chance constraints). Note that this methodology
has been also applied to the unit commitment problems \cite%
{wang2012chance,wu2014chance}. The classic solution to such stochastic
program is by scenario sampling \cite{Calafiore2005}, or robust optimization 
\cite{Mak2013}. In this paper, we adopt an approximation that is analogous
to the ``square-root staffing rule" from the queueing theory \cite{Whitt2007}%
, and derive a tractable approximation as mixed-integer second-order conic
program (MISOCP), which is less conservative and better interpreted than the
existing approaches.

\section{Model}

\label{s-model} \textbf{Energy network}. We consider an energy network
consisting of energy production facilities (\emph{plants}, denoted by a set $%
I$), to supply and satisfy local demands (\emph{customers}, denoted by a set 
$J$) of multiple types of energy. A plant may use different energy
production units (\emph{equipments}, denoted by a set $L$) to meet local
demands from customers. Customers generate demands for different \textit{%
types of energy} (denoted by a set $K$).

\textbf{Customers}. The customers' energy demands follow Poisson
distribution. We assume that orders from different demand locations are
independent. We parameterize the demand uncertainty by a \emph{scenario} $%
\omega $, which is drawn from the sample space $\Omega $ (collections of all
scenarios). Thus, we can write the demand from customer $j$ for energy $k$
under a given scenario $\omega $ as a Poisson random variable $\lambda
_{jk}\left( \omega \right) $. The expectation of a customer $j$'s demand for
energy $k$, which is the mean for the corresponding Poisson random variable,
is denoted by $\Lambda _{jk}$. Potential revenue from satisfying unit demand
for energy $k$ is $V_{k}$, regardless of customer.

\textbf{Plants}. A plant is to be placed at a chosen location, and the setup
cost for the plant at location $i$ is $f_{i}$. We use a binary decision
variable $Z_{i}$ to decide whether a plant should be in potential location $%
i $. We can generate the distance between customer $j$ and facility location 
$i $ as $d_{ij}$. The transportation cost $\phi _{k}\left( d_{ij}\right) $
can be any bounded function of distance $d_{ij}$, and also depending on the
energy type $k$. The cost of unit equipment $l$ to be invested at location $%
i $ is denoted by $g_{il}$, where the set of all equipments is $L$. $X_{il}$
is a decision variable representing the number of equipment $l$ available at
location $i$. Equipments can be placed at location $i$ only when a plant has
been set up at this location, i.e., 
\begin{equation*}
X_{il}\leq MZ_{i},\forall i\in I,l\in L,
\end{equation*}%
where $M$ is an arbitrarily large positive number, e.g., the maximum number
of equipments exhausting all investment.

\textbf{Multi-generation technology}. Since an equipment represents a
multi-generation unit, we use $\pi (k)$ to denote the set of equipments
which can be used to produce energy $k$, while $\pi ^{-1}(l)$ denotes the
types of energy which can be produced by equipment $l$. To model a \textit{%
flexible} energy production system, we assume a many-to-many mapping between
the types of energy and multi-generation units. Equipments differ in their
flexibility. The production capacity at each plant for energy $k$ is
determined by both the dedicated production units and the multi-generation
units which can be used to produce energy $k$.

\textbf{Production scheduling}. We define $Y_{ijk}(\omega )$ as a binary
variable if energy $k$ from customer $j$ should be produced from location $i$%
. Without loss of generality, the production rates of different equipments
are deterministic and are normalized to one. A feasible production schedule
should assign at least one plant to produce energy $k$ for customer $j$,
i.e., 
\begin{equation*}
\sum_{i\in I}Y_{ijk}\left( \omega \right) =1,\forall j\in J,k\in K,\omega
\in \Omega .
\end{equation*}

At the planning stage, the demands are not realized yet. The planner needs
to consider the uncertainty so that the equipments have sufficient capacity
for demands with high guarantee. Let $\xi \subseteq 2^{K}$ be any subset of
energy (wherein the notation $2^{K}$ represents the power set of $K$). The
aggregate production capacity for this set of energy is $\sum_{k\in \xi
}\sum_{l\in \pi (k)}X_{il}$, while the aggregate demand for this set of
energy is $\sum_{k\in \xi }\sum_{j}Y_{ijk}\left( \omega \right) \lambda
_{jk}\left( \omega \right) .$ For small $\alpha ,$ the \textit{service
guarantee constraints} (resource constraints) require that: 
\begin{equation*}
p\left\{ \sum_{k\in \xi }\sum_{j\in J}Y_{ijk}\left( \omega \right) \lambda
_{jk}\left( \omega \right) \leq \sum_{k\in \xi }\sum_{l\in \pi
(k)}X_{il}\right\} =1-\alpha ,\forall i\in I,\xi \subseteq 2^{K},
\end{equation*}%
which intuitively means that the aggregate capacity exceeds demand with high
probability. For given equipment investment decisions $\{X_{il}\}^{\prime }s$%
, these service guarantee constraints pin down feasible production schedule
decisions $Y_{ijk}\left( \omega \right) $.

Note that this formulation does not require the knowledge of which equipment
is producing energy for a particular customer. From a queueing theory
perspective, it is similar to a multi-server queueing model, which implies
that, a state-dependent dynamic scheduling policy will be required to \emph{%
implement} the optimization solution in practice \cite{bertsimas2010robust}.
We will provide further discussion on this issue in section \ref{s-comp}.

\textbf{Hierarchy of decision making}. Macro-level strategic decisions
include facility location and equipment investment. At the beginning, a
planner would decide how many plants are necessary, and where they should be
located. In addition, investment for different equipments are made.
Micro-level operational decisions include production and transportation
planning. After customers place orders, demands for energy are realized. The
planner decides which plant should produce energy to satisfy the demands.
Once production plans are scheduled, the corresponding transportation costs
are determined.

\textbf{Objective}. The planner maximizes the aggregate revenue subtracted
by the aggregate cost: 
\begin{equation*}
\mathbf{E}_{P}\sum_{i\in I,j\in J,k\in K}(V_{k}-\phi _{k}\left(
d_{ij}\right) )Y_{ijk}\left( \omega \right) \lambda _{jk}\left( \omega
\right) -\sum_{i\in I,l\in L}g_{il}X_{il}-\sum_{i\in I}f_{i}Z_{i}.
\end{equation*}%
Since the revenue is constant for each type of energy, it is equivalent to
adopt a cost-minimization framework, including transportation costs,
equipment investment costs, and plant setup costs, respectively. The
expectation is taken over probability $P$ associated with the sample space $%
\Omega $.

We summarize the nomenclature in this paper as follows.

\begin{supertabular}{p{1cm} p{6cm} p{1cm} p{6cm}  }
$\Omega $ & The set of all possible realization of uncertainty. &
$P$ & Probability measure for the uncertainty. \\

$I$ & The set of all potential locations for energy production. &
$J$ & The set of all customers (demand locations). \\

$K$ & The set of all types of energy. &
$L$ & The set of all equipments. \\

$V_{k}$ & Revenue from unit demand for energy $k$. &
$d_{ij}$ & Distance between customer $j$ and plant $i$.
\\

$\phi(\cdot) $ & Transportation cost function. &
$g_{il}$ & Cost of equipment $l$ at location $i$. \\

$f_{i}$ & Setup cost for a plant at potential location $i.$ &
$\lambda _{jk}\left( \cdot \right) $ & Demand at $j$ for energy $%
k $ (for given state or scenario). \\

$\Lambda _{jk}$ & Expected demand at $j$ for energy $k$. &
$Z_{i}$ & Binary decision variables for a plant at location $i$. \\

$Y_{ijk}(\cdot)$ & Whether energy $k$ at $j$ is to be
produced from $i$ (for given state or scenario). &
$X_{il}$ & Number (integer) of equipment $l$
available at location $i$. \\

$\pi (k)$ & The set of equipments to produce energy $k$. &
$\pi ^{-1}(l)$ & The types of energy produced by equipment $l$. \\
\end{supertabular}

\section{Analysis}

\label{s-analysis}In general, the problem is intractable due to the
uncountable set of scenarios $\Omega $. We now add more structure to the
probability space. Suppose that the demands are determined indirectly by
some finite underlying \textit{states}. In addition, we define $Y_{ijk}(s)$
as a binary decision variable whether the demand of energy $k$ from customer 
$j$ should be satisfied from the production at location $i$ when the state
is $s$. The set of all states is S and each state takes place independently
with probability $p_{s}$. Intuitively, a \textquotedblleft scenario" is a
random draw that can be viewed as \textquotedblleft raw data", while a $%
state $ represents market condition. The production decisions would be 
\textit{pooling} for each state. The motivation to use \textquotedblleft
state" in our formulation instead of \textquotedblleft scenario" is to
estimate \emph{average} demand intensity under different market conditions

\begin{lem}
Denote $\Lambda _{jk}\left( s\right) =E\left[ \lambda _{jk}\left( \omega
|s\right) \right] $, and assume that $\Lambda _{jk}\left( s\right)
\rightarrow \infty $. The resource constraint can be replaced by%
\begin{equation}
\Phi ^{-1}\left( 1-\alpha \right) \sqrt{\sum_{k\in \xi }\sum_{j\in
J}Y_{ijk}\left( s\right) \Lambda _{jk}\left( s\right) }+\sum_{k\in \xi
}\sum_{j\in J}Y_{ijk}\left( s\right) \Lambda _{jk}\left( s\right) \leq
\sum_{k\in \xi }\sum_{l\in \pi (k)}X_{il},\forall i\in I,s\in S,\xi
\subseteq 2^{K},  \tag{RE-$\xi $}
\end{equation}%
wherein $\Phi (\cdot )$ is the standard normal distribution function.
\end{lem}

Under such a resource investment policy, the system performance is analogous
to the \textit{quality-and-efficiency-driven regime} (QED) in queueing
theory \cite{Whitt2007}. Note that we need to assign at least one plant to
produce energy $k$ for customer $j$ for every state instead of every
scenario, i.e., we require that $\sum_{i\in I}Y_{ijk}\left( s\right)
=1,\forall j\in J,k\in K,s\in S.$ We summarize our model as follows:

\begin{equation*}
\text{Minimize }\sum_{i\in I,j\in J,k\in K,s\in S}\phi _{k}\left(
d_{ij}\right) Y_{ijk}\left( s\right) \Lambda _{jk}\left( s\right)
p_{s}+\sum_{i\in I,l\in L}g_{il}X_{il}+\sum_{i\in I}f_{i}Z_{i},
\end{equation*}%
subject to ($RE-\xi $), and%
\begin{equation}
\sum_{i\in I}Y_{ijk}\left( s\right) =1,\forall j\in J,k\in K,s\in S 
\tag{P-1}
\end{equation}%
\begin{equation*}
X_{il}\leq MZ_{i},\forall i\in I,l\in L
\end{equation*}%
\begin{equation*}
Y_{ijk}\left( s\right) ,Z_{i}\in \left\{ {0,1}\right\} ,X_{il}\in \mathbb{Z}.
\end{equation*}

\subsection{Flexibility Structures due to Multi-Generation Technologies}

In general, the capacity constraints RE-$\xi $ are needed for all possible
combinations of energy types. We need to eliminate the redundant constraints
to reduce the size of the optimization programs.

\begin{proposition}
Consider a partition $\xi _{t}\subset \xi ,$ for some $t=1,2,...$, $t<\infty 
$, wherein \{$\xi _{t}$\}'s are collectively exhaustive and mutually
exclusive. If $\cap _{t=1,2,..}\left( \cup _{k\in \xi _{t}}\pi (k)\right)
=\varnothing $, then RE-$\xi $ is redundant.
\end{proposition}

\textbf{Example 1} (Dense chaining). Suppose we have $\left\vert
L\right\vert $ types of equipments and $\left\vert K\right\vert =\left\vert
L\right\vert $ types of energy. We consider a particular dense flexibility
structure. Each equipment can be used to produce $\left\vert L\right\vert -1$
types of energy, i.e.,\ $\pi ^{-1}(l)=\{1,2,...,l-1,l+1,...,\left\vert
L\right\vert \}$, $\forall l=1,2,...,\left\vert L\right\vert -1.$ For $%
\forall i\in I,s\in S,$ the number of binding resource constraints can't be
greater than $\left\vert L\right\vert +1.$

\textbf{Example 2} (Star-flexibility). Suppose we have $\left\vert
L\right\vert $ types of equipment and $\left\vert K\right\vert =\left\vert
L\right\vert -1$ types of energy. Equipment $l=1,2,...,\left\vert
L\right\vert -1$ can only be used to produce one type of energy each, i.e., $%
\pi ^{-1}(\pi (k))=k,$ $\forall k=1,2,...,\left\vert L\right\vert -1.$
Equipment $l=\left\vert L\right\vert $ could be used to produce all types of
energy, i.e., $\pi ^{-1}(\left\vert L\right\vert )=K$. For $\forall i\in
I,s\in S,$ the number of binding resource constraints can't be greater than $%
2^{\left\vert L\right\vert -1}-1.$

\subsection{Capacity Allocation: Responsive vs. Anticipative}

\label{s-comp} We propose an alternative model in which the planner
explicitly allocates the production capacity for every equipments. Define $%
\Delta _{ilk}(s)$ as the proportion of time each equipment $l$ is scheduled
for energy $k$ at location $i.$ The model (P-1) can be modified as:

\begin{equation*}
\text{Minimize }\sum_{i\in I,j\in J,k\in K,s\in S}\phi _{k}\left(
d_{ij}\right) Y_{ijk}\left( s\right) \Lambda _{jk}\left( s\right)
p_{s}+\sum_{i\in I,l\in L}g_{il}X_{il}+\sum_{i\in I}f_{i}Z_{i},
\end{equation*}%
\begin{equation}
\sum_{i\in I}Y_{ijk}\left( s\right) =1,\forall j\in J,k\in K,s\in S 
\tag{P-2}
\end{equation}%
\begin{equation*}
X_{il}\leq MZ_{i},\forall i\in I,l\in L
\end{equation*}%
\begin{equation*}
\Phi ^{-1}\left( 1-\alpha \right) \sqrt{\sum_{j\in J}Y_{ijk}\left( s\right)
\Lambda _{jk}\left( s\right) }+\sum_{j\in J}Y_{ijk}\left( s\right) \Lambda
_{jk}\left( s\right) \leq \sum_{l\in \pi (k)}X_{il}\Delta _{ilk}(s),\forall
i\in I,k\in K,s\in S
\end{equation*}%
\begin{equation*}
\sum_{k}\Delta _{ilk}(s)\leq 1,\forall i\in I,l\in L,s\in S
\end{equation*}%
\begin{equation*}
Y_{ijk}\left( s\right) ,Z_{i}\in \left\{ {0,1}\right\} ,X_{il}\in \mathbb{Z}%
,\Delta _{ilk}(s)\in \lbrack 0,1].
\end{equation*}%
In what follows, we compare this model (P-2) with (P-1).

\begin{proposition}
The aggregate cost under \textit{anticipative allocation} (P-2) provides an
upper bound for that under \textit{responsive allocation} (P-1).
\end{proposition}

The names ``anticipative" and ``responsive" are motivated by the resemblance
of our resource constraints to the \emph{Conservation Law} in queueing
theory \citep{shanthikumar1992multiclass}. Using an analogy in queueing
systems, (P-2) is called ``anticipative" since we know how many servers
(equipments) will be used for certain types of customers (energy). (P-1) is
called ``responsive" since the optimization results only inform us what
equipments should be pooled together to produce what types of energy, the
implementation of which requires a dynamic scheduling mechanism.
Intuitively, with a \emph{fixed} production schedule in (P-2), the planner
fails to take advantage of the flexibility structure \emph{dynamically} (in 
\emph{response} to the realization of demands). Thus, the planner does not
make full use of the risk pooling instrument, and (P-2) yields a higher
aggregate cost.

\section{Computational Techniques and Numerical Examples}

\label{s-num} (P-1) is a mixed integer program, which is difficult to solve.
However, we can show that it can be solved easily under reformulation.

\begin{proposition}
(P-1) is equivalent to the following tractable MISOCP. 
\begin{equation*}
\text{Minimize}\sum_{i\in I,j\in J,k\in K,s\in S}\phi _{k}\left(
d_{ij}\right) Y_{ijk}\left( s\right) \Lambda _{jk}\left( s\right)
p_{s}+\sum_{i\in I,l\in L}g_{il}X_{il}+\sum_{i\in I}f_{i}Z_{i},
\end{equation*}%
\begin{equation}
\sum_{i\in I}Y_{ijk}\left( s\right) =1,\forall j\in J,k\in K,s\in S 
\tag{P-3}
\end{equation}%
\begin{equation*}
X_{il}\leq MZ_{i},\forall i\in I,l\in L
\end{equation*}%
\begin{equation*}
\sum_{k\in \xi }\sum_{j\in J}Y_{ijk}^{2}\left( s\right) \Lambda _{jk}\left(
s\right) \leqslant t_{i}^{2}\left( s\right) ,\forall i\in I,s\in S
\end{equation*}%
\begin{eqnarray*}
&&\left\Vert 
\begin{array}{c}
t_{i}\left( s\right) +\frac{\Phi ^{-1}\left( 1-\alpha \right) }{2} \\ 
\frac{1}{2}\left( 1-\sum_{k\in \xi }\sum_{l\in \pi (k)}X_{il}-\left[ \frac{%
\Phi ^{-1}\left( 1-\alpha \right) }{2}\right] ^{2}\right)%
\end{array}%
\right\Vert _{2} \\
&\leq &\frac{1}{2}\left( 1+\sum_{k\in \xi }\sum_{l\in \pi (k)}X_{il}+\left[ 
\frac{\Phi ^{-1}\left( 1-\alpha \right) }{2}\right] ^{2}\right) ,\forall
i\in I,s\in S,\xi \subseteq 2^{K}
\end{eqnarray*}%
\begin{equation*}
Y_{ijk}\left( s\right) ,Z_{i}\in \left\{ {0,1}\right\} ,t_{i}\left( s\right)
\in \mathbb{R}^{+},X_{il}\in \mathbb{Z}.
\end{equation*}
\end{proposition}

Using this computational technique, we conduct numerical experiments to
illustrate the economic value of the multi-generation technologies in lieu
of the spatial and inter-temporal demand uncertainty for energy: the average
demand is drawn from a two-type distribution wherein we interpret the demand
difference between two states as ``inter-temporal" while the difference
among locations as ``inter-spatial". We observe the following impacts of the
multi-generation technology from the numerical experiments:

\begin{itemize}
\item It reduces uncertainty by pooling capacity both within and across
different facilities.

\item It can compensate the loss in de-centralization due to the increase in
the transportation cost.

\item It can amplify the risk-pooling effects due to centralization
triggered by increasing setup costs.
\end{itemize}

The readers are referred to the Appendix for details of these major
observations summarized in the main paper.

\section{Conclusion}

In this paper, we investigate the system design problems in energy systems
when multi-generation technologies are available. Our modeling framework can
be adapted to incorporate extensions. For example, in many cases, we would
like to improve over the existing infrastructure instead of designing the
entire system. This case can be handled by replacing the corresponding
decision variables with known parameters. In addition, some operational
rules may have strategic impacts, such as some upper and lower limits for
production capacity and equipment start-up costs. These can also be taken
care of with additional operational constraints in the optimization problem.

In terms of future research, it will be fruitful to apply our methodology by
focusing on a particular sector, so that more structural results can be
generated to develop a problem-specific optimization algorithm. It will be
also desirable to implement our methodology in case studies for a specific
energy market. Furthermore, we assume in this paper that the demand is
uncertainty (state-dependent) while the energy production is stationary. For
the integration of certain sustainable energy, such as wind power, the
energy supply can also be uncertain (state-dependent). Another interesting
direction is to incorporate energy storage as an alternative instrument for
energy supply chain flexibility design. \label{s-con}

\bibliographystyle{chicago}
\bibliography{RML}

\begin{thebibliography}{}

\bibitem[\protect\citeauthoryear{Bertsimas and Doan}{Bertsimas and
  Doan}{2010}]{bertsimas2010robust}
Bertsimas, D. and X.~V. Doan (2010).
\newblock Robust and data-driven approaches to call centers.
\newblock {\em European Journal of Operational Research\/}~{\em 207\/}(2),
  1072--1085.

\bibitem[\protect\citeauthoryear{Calafiore and Campi}{Calafiore and
  Campi}{2005}]{Calafiore2005}
Calafiore, G. and M.~C. Campi (2005).
\newblock Uncertain convex programs: {R}andomized solutions and confidence
  levels.
\newblock {\em Mathematical Programming\/}~{\em 102\/}(1), 25--46.

\bibitem[\protect\citeauthoryear{Chen, Kang, O'Malley, Xia, Bai, Liu, Sun,
  Wang, and Li}{Chen et~al.}{2015}]{chen2015increasing}
Chen, X., C.~Kang, M.~O'Malley, Q.~Xia, J.~Bai, C.~Liu, R.~Sun, W.~Wang, and
  H.~Li (2015).
\newblock Increasing the flexibility of combined heat and power for wind power
  integration in {C}hina: Modeling and implications.
\newblock {\em {IEEE} Trans. Power Syst.\/}~{\em 30\/}(4), 1848--1857.

\bibitem[\protect\citeauthoryear{Chicco and Mancarella}{Chicco and
  Mancarella}{2009}]{chicco2009distributed}
Chicco, G. and P.~Mancarella (2009).
\newblock Distributed multi-generation: {A} comprehensive view.
\newblock {\em Renewable and Sustainable Energy Reviews\/}~{\em 13\/}(3),
  535--551.

\bibitem[\protect\citeauthoryear{Clegg and Mancarella}{Clegg and
  Mancarella}{2016}]{883bfc133a82486082cff0a5ea0e86fe}
Clegg, S. and P.~Mancarella (2016).
\newblock Integrated electrical and gas network flexibility assessment in
  low-carbon multi-energy systems.
\newblock {\em {IEEE} Trans. Sustain. Energy\/}~{\em 7\/}(2), 718--731.

\bibitem[\protect\citeauthoryear{Daskin}{Daskin}{2011}]{daskin2011network}
Daskin, M.~S. (2011).
\newblock {\em Network and discrete location: {M}odels, algorithms, and
  applications}.
\newblock John Wiley \& Sons.

\bibitem[\protect\citeauthoryear{Deng and Shen}{Deng and Shen}{2013}]{Deng2013}
Deng, T. and Z.-J.~M. Shen (2013).
\newblock Process flexibility design in unbalanced networks.
\newblock {\em Manufacturing \& Service Operations Management\/}~{\em 15\/}(1),
  24--32.

\bibitem[\protect\citeauthoryear{Li, Eremia, and Shahidehpour}{Li
  et~al.}{2008}]{li2008interdependency}
Li, T., M.~Eremia, and M.~Shahidehpour (2008).
\newblock Interdependency of natural gas network and power system security.
\newblock {\em {IEEE} Trans. Power Syst.\/}~{\em 23\/}(4), 1817--1824.

\bibitem[\protect\citeauthoryear{Ma, Kirschen, Belhomme, and Silva}{Ma
  et~al.}{2011}]{ma2011optimizing}
Ma, J., D.~S. Kirschen, R.~Belhomme, and V.~Silva (2011).
\newblock Optimizing the flexibility of a portfolio of generating plants.
\newblock In {\em Power Systems Computation Conference (PSCC)}.

\bibitem[\protect\citeauthoryear{Mak, Rong, and Shen}{Mak
  et~al.}{2013}]{Mak2013}
Mak, H.-Y., Y.~Rong, and Z.-J.~M. Shen (2013).
\newblock Infrastructure planning for electric vehicles with battery swapping.
\newblock {\em Management Science\/}~{\em 59\/}(7), 1557--1575.

\bibitem[\protect\citeauthoryear{Mak and Shen}{Mak and Shen}{2009a}]{Mak2009b}
Mak, H.-Y. and Z.-J.~M. Shen (2009a).
\newblock Stochastic programming approach to process flexibility design.
\newblock {\em Flexible Services and Manufacturing Journal\/}~{\em 21\/}(3-4),
  75--91.

\bibitem[\protect\citeauthoryear{Mak and Shen}{Mak and Shen}{2009b}]{Mak2009a}
Mak, H.-Y. and Z.-J.~M. Shen (2009b).
\newblock A two-echelon inventory-location problem with service considerations.
\newblock {\em Naval Research Logistics\/}~{\em 56\/}(8), 730--744.

\bibitem[\protect\citeauthoryear{Mancarella and Chicco}{Mancarella and
  Chicco}{2013}]{mancarella2013real}
Mancarella, P. and G.~Chicco (2013).
\newblock Real-time demand response from energy shifting in distributed
  multi-generation.
\newblock {\em {IEEE} Trans. Smart Grid\/}~{\em 4\/}(4), 1928--1938.

\bibitem[\protect\citeauthoryear{Martens}{Martens}{1998}]{martens1998energetic}
Martens, A. (1998).
\newblock The energetic feasibility of {CHP} compared to the separate
  production of heat and power.
\newblock {\em Applied Thermal Engineering\/}~{\em 18\/}(11), 935--946.

\bibitem[\protect\citeauthoryear{Martinez-Mares and
  Fuerte-Esquivel}{Martinez-Mares and
  Fuerte-Esquivel}{2012}]{martinez2012unified}
Martinez-Mares, A. and C.~R. Fuerte-Esquivel (2012).
\newblock A unified gas and power flow analysis in natural gas and electricity
  coupled networks.
\newblock {\em {IEEE} Trans. Power Syst.\/}~{\em 27\/}(4), 2156--2166.

\bibitem[\protect\citeauthoryear{Martinez-Mares and
  Fuerte-Esquivel}{Martinez-Mares and
  Fuerte-Esquivel}{2013}]{martinez2013robust}
Martinez-Mares, A. and C.~R. Fuerte-Esquivel (2013).
\newblock A robust optimization approach for the interdependency analysis of
  integrated energy systems considering wind power uncertainty.
\newblock {\em {IEEE} Trans. Power Syst.\/}~{\em 28\/}(4), 3964--3976.

\bibitem[\protect\citeauthoryear{Meunier}{Meunier}{2002}]{meunier2002co}
Meunier, F. (2002).
\newblock Co-and tri-generation contribution to climate change control.
\newblock {\em Applied Thermal Engineering\/}~{\em 22\/}(6), 703--718.

\bibitem[\protect\citeauthoryear{Qiu, Dong, Zhao, Meng, Zheng, and Hill}{Qiu
  et~al.}{2015}]{qiu2015low}
Qiu, J., Z.~Y. Dong, J.~H. Zhao, K.~Meng, Y.~Zheng, and D.~J. Hill (2015).
\newblock Low carbon oriented expansion planning of integrated gas and power
  systems.
\newblock {\em {IEEE} Trans. Power Syst.\/}~{\em 30\/}(2), 1035--1046.

\bibitem[\protect\citeauthoryear{Quelhas, Gil, McCalley, and Ryan}{Quelhas
  et~al.}{2007}]{quelhas2007multiperiod}
Quelhas, A., E.~Gil, J.~D. McCalley, and S.~M. Ryan (2007).
\newblock A multiperiod generalized network flow model of the {US} integrated
  energy system: Part {I} {M}odel description.
\newblock {\em {IEEE} Trans. Power Syst.\/}~{\em 22\/}(2), 829--836.

\bibitem[\protect\citeauthoryear{Shanthikumar and Yao}{Shanthikumar and
  Yao}{1992a}]{Shanthikumar1992}
Shanthikumar, J.~G. and D.~D. Yao (1992a).
\newblock Multiclass queueing systems: Polymatroidal structure and optimal
  scheduling control.
\newblock {\em Operations Research\/}~{\em 40\/}(3-supplement-2), S293--S299.

\bibitem[\protect\citeauthoryear{Shanthikumar and Yao}{Shanthikumar and
  Yao}{1992b}]{shanthikumar1992multiclass}
Shanthikumar, J.~G. and D.~D. Yao (1992b).
\newblock Multiclass queueing systems: Polymatroidal structure and optimal
  scheduling control.
\newblock {\em Operations Research\/}~{\em 40\/}(3-supplement-2), S293--S299.

\bibitem[\protect\citeauthoryear{Shen, Coullard, and Daskin}{Shen
  et~al.}{2003}]{Shen2003}
Shen, Z.-J.~M., C.~Coullard, and M.~S. Daskin (2003).
\newblock A joint location-inventory model.
\newblock {\em Transportation Science\/}~{\em 37\/}(1), 40--55.

\bibitem[\protect\citeauthoryear{Wang, Guan, and Wang}{Wang
  et~al.}{2012}]{wang2012chance}
Wang, Q., Y.~Guan, and J.~Wang (2012).
\newblock A chance-constrained two-stage stochastic program for unit commitment
  with uncertain wind power output.
\newblock {\em {IEEE} Trans. Power Syst.\/}~{\em 27\/}(1), 206--215.

\bibitem[\protect\citeauthoryear{Whitt}{Whitt}{2007}]{Whitt2007}
Whitt, W. (2007).
\newblock What you should know about queueing models to set staffing
  requirements in service systems.
\newblock {\em Naval Research Logistics (NRL)\/}~{\em 54\/}(5), 476--484.

\bibitem[\protect\citeauthoryear{Wu, Shahidehpour, Li, and Tian}{Wu
  et~al.}{2014}]{wu2014chance}
Wu, H., M.~Shahidehpour, Z.~Li, and W.~Tian (2014).
\newblock Chance-constrained day-ahead scheduling in stochastic power system
  operation.
\newblock {\em {IEEE} Trans. Power Syst.\/}~{\em 29\/}(4), 1583--1591.

\bibitem[\protect\citeauthoryear{Zipkin}{Zipkin}{1980}]{Zipkin1980}
Zipkin, P.~H. (1980).
\newblock Bounds for aggregating nodes in network problems.
\newblock {\em Mathematical Programming\/}~{\em 19\/}(1), 155--177.

\end{thebibliography}

\appendix

\section{Appendix. Proofs.}

In this appendix, we first provide the detailed proofs of the main
results.

\textbf{Proof of Lemma 1. } With the scenario aggregation that $\Lambda
_{jk}\left( s\right) =E\left[ \lambda _{jk}\left( \omega |s\right) \right] $%
, we can re-write the resource constraint as 
\begin{equation}
p\left\{ \sum_{k\in \xi }\sum_{j}Y_{ijk}\left( s\right) \lambda _{jk}\left(
\omega |s\right) \leq \sum_{k\in \xi }\sum_{l\in \pi (k)}X_{il}\right\}
=1-\alpha ,\forall s\in S.
\end{equation}%
Next we show that this constraint can be replaced by

\begin{equation}
\Phi ^{-1}\left( 1-\alpha \right) \sqrt{\sum_{k\in \xi }\sum_{j\in
J}Y_{ijk}\left( s\right) \Lambda _{jk}\left( s\right) }+\sum_{k\in \xi
}\sum_{j\in J}Y_{ijk}\left( s\right) \Lambda _{jk}\left( s\right) \leq
\sum_{k\in \xi }\sum_{l\in \pi (k)}X_{il},\forall i\in I,s\in S,\xi
\subseteq 2^{K},  \tag{RE-$\xi $}
\end{equation}%
as $\Lambda _{jk}\left( s\right) \rightarrow \infty $, wherein $\Phi $ is
the distribution function for standard normal random variable. To see this,
we replace the service guarantee by Gaussian approximation:

\begin{equation}
p\left\{ \sum_{k\in \xi }\sum_{j\in J}Y_{ijk}\left( s\right) \lambda
_{jk}\left( \omega |s\right) \leq \sum_{k\in \xi }\sum_{l\in \pi
(k)}X_{il}\right\} =p\left\{ 
\begin{array}{c}
\frac{\sum_{k\in \xi }\sum_{j\in J}Y_{ijk}\left( s\right) [\lambda
_{jk}\left( \omega |s\right) -\Lambda _{jk}\left( s\right) ]}{\sqrt{%
\sum_{k\in \xi }\sum_{j\in J}Y_{ijk}\left( s\right) \Lambda _{jk}\left(
s\right) }} \\ 
\leq \frac{\sum_{k\in \xi }\sum_{l\in \pi (k)}X_{il}-\sum_{k\in \xi
}\sum_{j\in J}Y_{ijk}\left( s\right) \Lambda _{jk}\left( s\right) }{\sqrt{%
\sum_{k\in \xi }\sum_{j\in J}Y_{ijk}\left( s\right) \Lambda _{jk}\left(
s\right) }},%
\end{array}%
\right\}
\end{equation}%
for $\forall i\in I,s\in S,\xi \subseteq 2^{K}.$ Let a feasible investment
level be to choose the minimum $\sum_{k\in \xi }\sum_{l\in \pi (k)}X_{il}$
that is larger than $\sum_{k\in \xi }\sum_{j}Y_{ijk}\left( s\right) \Lambda
_{jk}\left( s\right) +\Phi ^{-1}\left( 1-\alpha \right) \sqrt{\sum_{k\in \xi
}\sum_{j}Y_{ijk}\left( s\right) \Lambda _{jk}\left( s\right) }$, the
right-hand side is $1-\alpha $ due to the Central Limit Theorem. $\square $

\textbf{Proof of Proposition 1.} For $\forall \xi ,$ $\sum_{k\in \xi
}\sum_{j}Y_{ijk}\left( s\right) \Lambda _{jk}\left( s\right) +\Phi
^{-1}\left( 1-\alpha \right) \sqrt{\sum_{k\in \xi }\sum_{j}Y_{ijk}\left(
s\right) \Lambda _{jk}\left( s\right) }\newline
$

\begin{eqnarray}
&=&\left( \sum_{k\in \xi _{1}}+\sum_{k\in \xi _{2}}+...\right)
\sum_{j}Y_{ijk}\left( s\right) \Lambda _{jk}\left( s\right) +\Phi
^{-1}\left( 1-\alpha \right) \sqrt{\left( \sum_{k\in \xi _{1}}+\sum_{k\in
\xi _{2}}+...\right) \sum_{j}Y_{ijk}\left( s\right) \Lambda _{jk}\left(
s\right) }\newline
\notag \\
&\leq &\sum_{k\in \xi }\sum_{l\in \pi (k)}X_{il},\forall i\in I,s\in S.
\end{eqnarray}%
Therefore, RE-$\xi $ is redundant. Note that the converse is not true, i.e.,
if $\cap _{k\in \xi }\pi (k)\neq \varnothing $, then it is still possible
that RE-$\xi $ is redundant. $\square $

\textbf{Proof of Proposition 2.} To prove that this is a upper bound, it
suffices to check the feasibility of (P-2)'s solution to (P-1)'s constraints:%
\begin{eqnarray}
&&\sum_{k\in \xi }\sum_{j}Y_{ijk}\left( s\right) \Lambda _{jk}\left(
s\right) +\Phi ^{-1}\left( 1-\alpha \right) \sqrt{\sum_{k\in \xi
}\sum_{j}Y_{ijk}\left( s\right) \Lambda _{jk}\left( s\right) }  \notag \\
&\leq &\sum_{k\in \xi }\sum_{j}Y_{ijk}\left( s\right) \Lambda _{jk}\left(
s\right) +\Phi ^{-1}\left( 1-\alpha \right) \sum_{k\in \xi }\sqrt{%
\sum_{j}Y_{ijk}\left( s\right) \Lambda _{jk}\left( s\right) }  \notag \\
&\leq &\sum_{k\in \xi }\sum_{l\in \pi (k)}X_{il}\Delta _{ilk}(s)  \notag \\
&\leq &\sum_{k\in \xi }\sum_{l\in \pi (k)}X_{il},\forall i\in I,\xi \in
2^{K},s\in S.
\end{eqnarray}%
The first inequality is due to the property of square-root function (or
Jensen's inequality). The second inequality holds via the feasibility of
(P-2). The third inequality is true since $\Delta _{ilk}(s)\leq 1$. Since
any solution of (P-2) is feasible for (P-1), (P-2) is more restrictive and
thus returns a higher cost. $\square $

\textbf{Proof of Proposition 3.} For convenience, we write (P-1) as follows: 
\begin{equation*}
\text{Minimize }\sum_{i\in I,j\in J,k\in K,s\in S}\phi _{k}\left(
d_{ij}\right) Y_{ijk}\left( s\right) \Lambda _{jk}\left( s\right)
p_{s}+\sum_{i\in I,l\in L}g_{il}X_{il}+\sum_{i\in I}f_{i}Z_{i},
\end{equation*}%
subject to 
\begin{equation}
\Phi ^{-1}\left( 1-\alpha \right) \sqrt{\sum_{k\in \xi }\sum_{j\in
J}Y_{ijk}\left( s\right) \Lambda _{jk}\left( s\right) }+\sum_{k\in \xi
}\sum_{j\in J}Y_{ijk}\left( s\right) \Lambda _{jk}\left( s\right) \leq
\sum_{k\in \xi }\sum_{l\in \pi (k)}X_{il},\forall i\in I,s\in S,\xi
\subseteq 2^{K}.  \notag
\end{equation}%
\begin{equation}
\sum_{i\in I}Y_{ijk}\left( s\right) =1,\forall j\in J,k\in K,s\in S 
\tag{P-1}
\end{equation}%
\begin{equation*}
X_{il}\leq MZ_{i},\forall i\in I,l\in L
\end{equation*}%
\begin{equation*}
Y_{ijk}\left( s\right) ,Z_{i}\in \left\{ {0,1}\right\} ,X_{il}\in \mathbf{Z}.
\end{equation*}

Since $Y_{ijk}\left( s\right) ^{\prime }s$ are binary, they will be
equivalent to $Y_{ijk}^{2}\left( s\right) $, by which we can re-write the
resource constraints as

\begin{equation}
\Phi ^{-1}\left( 1-\alpha \right) \sqrt{\sum_{k\in \xi }\sum_{j\in
J}Y_{ijk}^{2}\left( s\right) \Lambda _{jk}\left( s\right) }+\sum_{k\in \xi
}\sum_{j\in J}Y_{ijk}^{2}\left( s\right) \Lambda _{jk}\left( s\right) \leq
\sum_{k\in \xi }\sum_{l\in \pi (k)}X_{il},\forall i\in I,s\in S,\xi
\subseteq 2^{K}.
\end{equation}

We can introduce variable $t_{i}\left( s\right) ^{\prime }s$ such that 
\begin{equation}
\sum_{k\in \xi }\sum_{j\in J}Y_{ijk}^{2}\left( s\right) \Lambda _{jk}\left(
s\right) =t_{i}^{2}\left( s\right) ,
\end{equation}%
we can re-write the resource constraints as

\begin{equation}
\Phi ^{-1}\left( 1-\alpha \right) t_{i}\left( s\right) +t_{i}^{2}\left(
s\right) \leq \sum_{k\in \xi }\sum_{l\in \pi (k)}X_{il},\forall i\in I,s\in
S,\xi \subseteq 2^{K}.
\end{equation}%
Adding $\left[ \frac{\Phi ^{-1}\left( 1-\alpha \right) }{2}\right] ^{2}$ on
both sides:

\begin{equation}
\left[ t_{i}\left( s\right) +\frac{\Phi ^{-1}\left( 1-\alpha \right) }{2}%
\right] ^{2}\leq \sum_{k\in \xi }\sum_{l\in \pi (k)}X_{il}+\left[ \frac{\Phi
^{-1}\left( 1-\alpha \right) }{2}\right] ^{2},\forall i\in I,s\in S,\xi
\subseteq 2^{K}.
\end{equation}%
Now the problem becomes:

\begin{equation*}
\text{Minimize }\sum_{i\in I,j\in J,k\in K,s\in S}\phi _{k}\left(
d_{ij}\right) Y_{ijk}\left( s\right) \Lambda _{jk}\left( s\right)
p_{s}+\sum_{i\in I,l\in L}g_{il}X_{il}+\sum_{i\in I}f_{i}Z_{i},
\end{equation*}%
subject to

\begin{equation}
\sum_{i\in I}Y_{ijk}\left( s\right) =1,\forall j\in J,k\in K,s\in S 
\tag{P-1'}
\end{equation}%
\begin{equation*}
X_{il}\leq MZ_{i},\forall i\in I,l\in L
\end{equation*}%
\begin{equation*}
\sum_{k\in \xi }\sum_{j\in J}Y_{ijk}^{2}\left( s\right) \Lambda _{jk}\left(
s\right) =t_{i}^{2}\left( s\right) ,\forall i\in I,s\in S
\end{equation*}%
\begin{equation*}
\left[ t_{i}\left( s\right) +\frac{\Phi ^{-1}\left( 1-\alpha \right) }{2}%
\right] ^{2}\leq \sum_{k\in \xi }\sum_{l\in \pi (k)}X_{il}+\left[ \frac{\Phi
^{-1}\left( 1-\alpha \right) }{2}\right] ^{2},\forall i\in I,s\in S,\xi
\subseteq 2^{K}.
\end{equation*}%
\begin{equation*}
Y_{ijk}\left( s\right) ,Z_{i}\in \left\{ {0,1}\right\} ,t_{i}\left( s\right)
\in \mathbb{R}^{+},X_{il}\in \mathbb{Z}.
\end{equation*}%
We shall now prove that \textit{(P-1') and (P-3) }are equivalent by
contradiction\textit{. }We write down \textit{(P-3) }for your convenience.

\begin{equation*}
\text{Minimize}\sum_{i\in I,j\in J,k\in K,s\in S}\phi _{k}\left(
d_{ij}\right) Y_{ijk}\left( s\right) \Lambda _{jk}\left( s\right)
p_{s}+\sum_{i\in I,l\in L}g_{il}X_{il}+\sum_{i\in I}f_{i}Z_{i},
\end{equation*}%
\begin{equation}
\sum_{i\in I}Y_{ijk}\left( s\right) =1,\forall j\in J,k\in K,s\in S 
\tag{P-3}
\end{equation}%
\begin{equation*}
X_{il}\leq MZ_{i},\forall i\in I,l\in L
\end{equation*}%
\begin{equation*}
\sum_{j\in J}Y_{ijk}^{2}\left( s\right) \Lambda _{jk}\left( s\right)
\leqslant t_{i}^{2}\left( s\right) ,\forall i\in I,s\in S
\end{equation*}%
\begin{equation*}
\left[ t_{i}\left( s\right) +\frac{\Phi ^{-1}\left( 1-\alpha \right) }{2}%
\right] ^{2}\leq \sum_{k\in \xi }\sum_{l\in \pi (k)}X_{il}+\left[ \frac{\Phi
^{-1}\left( 1-\alpha \right) }{2}\right] ^{2},\forall i\in I,s\in S,\xi
\subseteq 2^{K},
\end{equation*}%
\begin{equation*}
Y_{ijk}\left( s\right) ,Z_{i}\in \left\{ {0,1}\right\} ,t_{i}\left( s\right)
\in \mathbb{R}^{+},X_{il}\in \mathbb{Z}.
\end{equation*}%
Firstly, any solution of \textit{(P-1') satisfies (P-3)}, since \textit{(P-3)%
} is a relaxation. Suppose that a solution of\textit{\ (P-3)} does not
satisfy \textit{(P-1')}, then $\sum_{k\in \xi }\sum_{j\in
J}Y_{ijk}^{2}\left( s\right) \Lambda _{jk}\left( s\right) =t_{i}^{2}\left(
s\right) ,\forall i\in I,s\in S$, is not strictly binding: $\exists i\in
I,s\in S$ such that $\sum_{k\in \xi }\sum_{j\in J}Y_{ijk}^{2}\left( s\right)
\Lambda _{jk}\left( s\right) <t_{i}^{2}\left( s\right) $. Then, suppose that 
$t_{i}^{2}\left( s\right) =\sum_{k\in \xi }\sum_{j\in J}Y_{ijk}^{2}\left(
s\right) \Lambda _{jk}\left( s\right) +\varepsilon $ for some $\varepsilon
>0 $. However, this can not be optimal for \textit{(P-3)}, since we can
further decrease $t_{i}\left( s\right) $. As

\begin{equation}
\left[ t_{i}\left( s\right) +\frac{\Phi ^{-1}\left( 1-\alpha \right) }{2}%
\right] ^{2}\leq \sum_{k\in \xi }\sum_{l\in \pi (k)}X_{il}+\left[ \frac{\Phi
^{-1}\left( 1-\alpha \right) }{2}\right] ^{2},
\end{equation}%
we can reduce some of the $X_{il}$ and thus the objective function. \textit{%
Reductio ad absurdum}, solution of\textit{\ (P-3)} satisfy \textit{(P-1')},
and thus they are equivalent.

Finally, we need to show that

\begin{equation}
\Phi ^{-1}\left( 1-\alpha \right) t_{i}\left( s\right) +t_{i}^{2}\left(
s\right) \leq \sum_{k\in \xi }\sum_{l\in \pi (k)}X_{il},\forall i\in I,s\in
S,\xi \subseteq 2^{K},
\end{equation}%
is a second-order cone. We have already demonstrated that it is equivalent to

\begin{equation}
\left[ t_{i}\left( s\right) +\frac{\Phi ^{-1}\left( 1-\alpha \right) }{2}%
\right] ^{2}\leq \sum_{k\in \xi }\sum_{l\in \pi (k)}X_{il}+\left[ \frac{\Phi
^{-1}\left( 1-\alpha \right) }{2}\right] ^{2},\forall i\in I,s\in S,\xi
\subseteq 2^{K},
\end{equation}%
which is quadratic, and can be re-written as:

\begin{equation}
\left\Vert 
\begin{array}{c}
t_{i}\left( s\right) +\frac{\Phi ^{-1}\left( 1-\alpha \right) }{2} \\ 
\frac{1}{2}\left( 1-\sum_{k\in \xi }\sum_{l\in \pi (k)}X_{il}-\left[ \frac{%
\Phi ^{-1}\left( 1-\alpha \right) }{2}\right] ^{2}\right)%
\end{array}%
\right\Vert _{2}\leq \frac{1}{2}\left( 1+\sum_{k\in \xi }\sum_{l\in \pi
(k)}X_{il}+\left[ \frac{\Phi ^{-1}\left( 1-\alpha \right) }{2}\right]
^{2}\right) .
\end{equation}%
$\square $

\section{Numerical Results.}

We assume that there are three different types of energy, i.e., $K=\{A,S,H\}$%
. There are four different equipments, i.e., $L=\{A,S,H,F\}$, where $F$
indicates an equipment with full flexibility. The plants are located
potentially in $J$, a set of 18 cities (LA, SD, SJ, SF, LB, OAK, SAC, FRE,
SA, ANA, RIV, STKN, HB, GNDL, BKD, FMT, MOD, SB). Demands are also from the
same potential locations, i.e., $I=J$. There are two possible states of the
world, i.e., $S=\{H,L\}$.

Transportation costs are generated proportional to the distance $d_{ij}$
from \cite{daskin2011network}, i.e., $\phi _{k}\left( d_{ij}\right) =\phi
d_{ij}$ and the scaling coefficient $\phi =2$. We assume that the investment
coefficient $\Phi ^{-1}\left( 1-\alpha \right) =0.2$. We use heterogeneous
service rates as a generalization of the model, wherein (P1) can be easily
adapted by the following modification to RE-$\xi $ constraints:%
\begin{equation}
\Phi ^{-1}\left( 1-\alpha \right) \sqrt{\sum_{k\in \xi }\sum_{j\in
J}Y_{ijk}\left( s\right) \Lambda _{jk}\left( s\right) }+\sum_{k\in \xi
}\sum_{j\in J}Y_{ijk}\left( s\right) \Lambda _{jk}\left( s\right) \leq
\sum_{k\in \xi }\sum_{l\in \pi (k)}\mu _{ilk}X_{il},\forall i\in I,s\in
S,\xi \subseteq 2^{K},  \notag
\end{equation}%
wherein the service rates can depend on location $i$, the equipment $l$, and
potentially the energy type $k$, without changing the structure of this
optimization problem. We assume that $\mu _{A}=500$, $\mu _{S}=460$, $\mu
_{H}=420$. The multi-generation technologies can be used to produce all
three types of energy, and the corresponding service rate is $\mu _{F}=400$.
Two states are realized with equal probability, i.e., $\Pr (H)=\Pr (L)=0.5$.
The setup costs for potential locations are proportional to the real estate
prices in that city.

We consider a revenue maximizing framework, wherein valuations for energy
consumption are heterogeneous, i.e., $V_{A}=800$, $V_{S}=725$, $V_{H}=600.$
Equipment costs are homogeneous across all locations, i.e., $g_{iA}=90000$, $%
g_{iS}=70000$, $g_{iH}=50000$, $g_{iF}=100000$, $\forall i\in I.$ Finally,
we generate aggregate demand in proportion to the population in each city,
and split across different types of energy. We choose to split in a way to
generate both spatial and inter-temporal (state-dependent) heterogeneity
(see appendix for the demand distribution). 
\begin{table}[htbp]
\caption{Optimal equipments quantities, with multi-generation technologies.}
\label{tab:spatialheterflex}\centering
\begin{tabular}{cccccccccc}
\toprule & \multicolumn{4}{c}{High spatial heterogeneity} &  & 
\multicolumn{4}{c}{Low spatial heterogeneity} \\ 
\midrule Location Chosen & LA & MOD & OAK & SD &  & FRE & LA & OAK & SD \\ 
\# Equipment H & 5 & 1 & 2 & 2 &  & 0 & 5 & 3 & 2 \\ 
\# Equipment A & 4 & 0 & 2 & 1 &  & 0 & 4 & 2 & 1 \\ 
\# Equipment S & 0 & 0 & 0 & 0 &  & 0 & 2 & 1 & 1 \\ 
\# Multi-generation units & 3 & 1 & 2 & 1 &  & 1 & 1 & 1 & 0 \\ 
\bottomrule &  &  &  &  &  &  &  &  & 
\end{tabular}%
\end{table}

\begin{table}[htbp]
\caption{Optimal equipments quantities, without multi-generation
technologies.}
\label{tab:spatialheternonflex}\centering
\begin{tabular}{cccccccccc}
\toprule & \multicolumn{4}{c}{High spatial heterogeneity} &  & 
\multicolumn{4}{c}{Low spatial heterogeneity} \\ 
\midrule Location Chosen & FRE & LA & OAK & SD &  & LA & MOD & OAK & SD \\ 
\# Equipment H & 1 & 6 & 4 & 3 &  & 5 & 1 & 3 & 2 \\ 
\# Equipment A & 1 & 6 & 4 & 2 &  & 5 & 1 & 2 & 1 \\ 
\# Equipment S & 1 & 2 & 2 & 1 &  & 2 & 1 & 1 & 1 \\ 
\bottomrule &  &  &  &  &  &  &  &  & 
\end{tabular}%
\end{table}

Table \ref{tab:spatialheterflex} and Table \ref{tab:spatialheternonflex}
summarize the optimal equipment quantities for the chosen locations under
both high and low spatial demand heterogeneity. We compare the results with
and without multi-generation technologies. Under high spatial heterogeneity,
it would be optimal to invest in the multi-generation technologies more than
that under low spatial heterogeneity. In addition, with the investment in
multi-generation technologies, the number of equipments needed is smaller
than that in the non-flexible case. Such improvement is more significant
when spatial heterogeneity is high. We summarize the insight as follows.

\begin{observation}
The multi-generation technologies reduce uncertainty by pooling capacity
both within and across different facilities.
\end{observation}

\begin{table}[htbp]
\caption{Sensitivity analysis for transportation cost coefficient, with
multi-generation technologies}
\label{tab:transflex}\centering
\begin{tabular}{rrrrrr}
\toprule $\phi$ & 1 & 1.5 & 2 & 2.5 & 3 \\ 
\midrule Total \# open sites & 3 & 5 & 7 & 7 & 7 \\ 
Total \# multi-generation units & 2 & 2 & 5 & 5 & 5 \\ 
Total \# equipments H & 9 & 9 & 9 & 9 & 9 \\ 
Total \# equipments A & 15 & 15 & 14 & 14 & 14 \\ 
Total \# equipment S & 4 & 4 & 3 & 3 & 3 \\ 
Total equipment cost (\$) & 2280000 & 2280000 & 2420000 & 2420000 & 2420000
\\ 
Total transportation cost (\$) & 437045 & 452159 & 335854 & 378496 & 454195
\\ 
Total net revenue (\$) & 4648697 & 4483583 & 4359888 & 4277246 & 4201547 \\ 
\bottomrule &  &  &  &  & 
\end{tabular}%
\end{table}

\begin{table}[htbp]
\caption{Sensitivity analysis for transportation cost coefficient, without
multi-generation technologies}
\label{tab:transnonflex}\centering
\begin{tabular}{rrrrrr}
\toprule $\phi$ & 1 & 1.5 & 2 & 2.5 & 3 \\ 
\midrule Total \# open sites & 3 & 4 & 6 & 6 & 6 \\ 
Total \# equipments H & 9 & 10 & 11 & 11 & 11 \\ 
Total \# equipments A & 17 & 17 & 17 & 17 & 17 \\ 
Total \# equipment S & 5 & 5 & 6 & 6 & 6 \\ 
Aggregate equipment cost (\$) & 2330000 & 2380000 & 2500000 & 2500000 & 
2500000 \\ 
Aggregate transportation cost (\$) & 469570 & 582149 & 498881 & 623601 & 
748321 \\ 
Aggregate net revenue (\$) & 4516172 & 4303593 & 4166861 & 4042141 & 3917421
\\ 
\bottomrule &  &  &  &  & 
\end{tabular}%
\end{table}
Table \ref{tab:transflex} and Table \ref{tab:transnonflex} compare the
optimal facility location and equipment investment decisions with different
transportation costs, with and without the multi-generation technologies.
With the increase of the transportation cost coefficient $\phi$, it is
optimal to set up more plants so that the demands can be satisfied by the
nearest location. With the increasing number of the open sites, more
multi-generation units are needed because the production capacity is more
de-centralized, and the system suffers more from demand uncertainty both
within and across different locations. Surprisingly, the aggregate
transportation cost decreases when $\phi$ increases from $1.5$ to $2$.

\begin{observation}
Investment in the multi-generation technologies can compensate the loss in
de-centralization due to the increase of the transportation cost
coefficient. Consequently, the overall transportation costs could decrease.
\end{observation}

The increase of the setup cost has the opposite effect to the transportation
cost coefficient, summarized in the following observation and the remaining
two tables. 
\begin{table}[htbp]
\caption{Sensitivity analysis for set up cost, with multi-generation
technologies.}
\label{tab:setupsenseflex}\centering
\begin{tabular}{cccccc}
\toprule Set up cost & Very low & Low & Medium & High & Very high \\ 
\midrule Total \# open sites & 7 & 4 & 3 & 3 & 2 \\ 
Total \# multi-generation units & 5 & 5 & 2 & 2 & 3 \\ 
Total \# equipments H & 9 & 11 & 9 & 9 & 10 \\ 
Total \# equipments A & 14 & 13 & 15 & 15 & 14 \\ 
Total \# equipment S & 3 & 2 & 4 & 4 & 3 \\ 
Aggregate equipment cost (\$) & 2420000 & 2360000 & 2280000 & 2280000 & 
2270000 \\ 
Aggregate transportation cost (\$) & 335854 & 605403 & 776704 & 925164 & 
1377640 \\ 
Aggregate net revenue (\$) & 4359888 & 3931859 & 3677198 & 3482578 & 3095602
\\ 
\bottomrule &  &  &  &  & 
\end{tabular}%
\end{table}

\begin{table}[htbp]
\caption{Sensitivity analysis for set up cost, without multi-generation
technologies.}
\label{tab:setupnonflex}\centering
\begin{tabular}{cccccc}
\toprule Set up cost & Very low & Low & Medium & High & Very high \\ 
\midrule Total \# open sites & 6 & 4 & 4 & 3 & 2 \\ 
Total \# equipments H & 11 & 10 & 10 & 10 & 10 \\ 
Total \# equipments A & 17 & 17 & 17 & 17 & 17 \\ 
Total \# equipment S & 6 & 6 & 6 & 6 & 6 \\ 
Aggregate equipment cost (\$) & 2500000 & 2450000 & 2450000 & 2450000 & 
2450000 \\ 
Aggregate transportation cost (\$) & 498881 & 683499 & 683499 & 945895 & 
1366350 \\ 
Aggregate net revenue (\$) & 4166861 & 3771763 & 3501403 & 3291847 & 2926894
\\ 
\bottomrule &  &  &  &  & 
\end{tabular}%
\end{table}

\begin{observation}
The multi-generation technologies amplify the risk-pooling effects due to
centralization triggered by the increase of the setup costs.
\end{observation}

\section{Demand Data.}

In this appendix, we provide the demand data we use in the numerical
experiment. The demand is drawn from a two-type distribution wherein we
interpret the demand difference between states ``H" and ``L" as
``inter-temporal" while the difference among locations as ``inter-spatial".
We also include the location-dependent setup costs used in our calculations. 
\begin{table}[tbph]
\caption{Setup costs for each potential locations.}
\label{tab:setupcost}\centering
\begin{tabular}{ccccccc}
\toprule City & LA & SD & SJ & SF & LB & OAK \\ 
Setup Cost (\$) & 214880 & 178320 & 238000 & 328400 & 151600 & 159280 \\ 
\hline
City & SAC & FRE & SA & ANA & RIV & STKN \\ 
Setup Cost (\$) & 106800 & 86000 & 152720 & 165040 & 113200 & 78000 \\ \hline
City & HB & GNDL & BKD & FMT & MOD & SB \\ 
Setup Cost (\$) & 233200 & 202800 & 93040 & 230000 & 85040 & 79280 \\ 
\bottomrule &  &  &  &  &  & 
\end{tabular}%
\end{table}
\begin{table}[htbp]
\caption{Demands under high spatio-temporal heterogeneity.}
\label{tab:demandhighheter}\centering
\begin{tabular}{crrrccc}
\toprule States & \multicolumn{3}{c}{H} & \multicolumn{3}{c}{L} \\ 
\midrule Cities & \multicolumn{1}{c}{Demands (A)} & \multicolumn{1}{c}{
Demands (S)} & \multicolumn{1}{c}{Demands (H)} & Demands (A) & Demands (S) & 
Demands (H) \\ 
LA & 189.19 & 26.2764 & 47.2975 & 42.0422 & 157.658 & 63.0634 \\ 
SD & 134.929 & 18.7401 & 33.7322 & 18.7401 & 93.7005 & 74.9604 \\ 
SF & 43.3426 & 126.416 & 10.8356 & 108.356 & 18.0594 & 54.1782 \\ 
SJ & 263.103 & 36.5421 & 65.7758 & 36.5421 & 182.71 & 146.168 \\ 
LB & 148.185 & 20.5812 & 37.0462 & 20.5812 & 102.906 & 82.3248 \\ 
OAK & 24.7964 & 144.646 & 37.1947 & 144.646 & 20.6637 & 41.3274 \\ 
SAC & 438.266 & 2556.55 & 657.399 & 2191.33 & 730.443 & 730.443 \\ 
FRE & 343.386 & 47.6925 & 85.8465 & 76.308 & 286.155 & 114.462 \\ 
ANA & 41.7298 & 121.712 & 10.4324 & 104.324 & 17.3874 & 52.1622 \\ 
RIV & 104.055 & 303.494 & 26.0138 & 260.138 & 43.3563 & 130.069 \\ 
SA & 163 & 22.6389 & 40.75 & 36.2222 & 135.833 & 54.3334 \\ 
HB & 25.78 & 150.383 & 38.6699 & 150.383 & 21.4833 & 42.9666 \\ 
STKN & 104 & 303.332 & 25.9999 & 259.999 & 43.3332 & 130 \\ 
GNDL & 117.681 & 16.3446 & 29.4203 & 16.3446 & 81.723 & 65.3784 \\ 
BKD & 146.195 & 852.802 & 219.292 & 852.802 & 121.829 & 243.658 \\ 
FMT & 220.02 & 641.726 & 55.0051 & 550.051 & 91.6752 & 275.026 \\ 
MOD & 180.157 & 525.458 & 45.0392 & 450.392 & 75.0654 & 225.196 \\ 
SB & 49.5317 & 144.467 & 12.3829 & 123.829 & 20.6382 & 61.9146 \\ 
\bottomrule &  &  &  &  &  & 
\end{tabular}%
\end{table}

\begin{table}[htbp]
\caption{Demands under low spatio-temporal heterogeneity.}
\label{tab:demandlowheter}\centering
\begin{tabular}{crrr}
\toprule States & \multicolumn{3}{c}{Both H and L} \\ 
\midrule Cities & \multicolumn{1}{c}{Demands (A)} & \multicolumn{1}{c}{
Demands (S)} & \multicolumn{1}{c}{Demands (H)} \\ 
LA & 115.616 & 91.9674 & 55.1804 \\ 
SD & 76.8344 & 56.2203 & 54.3463 \\ 
SF & 75.8495 & 72.2376 & 32.5069 \\ 
SJ & 149.823 & 109.626 & 105.972 \\ 
LB & 84.3829 & 61.7436 & 59.6855 \\ 
OAK & 84.7212 & 82.6548 & 39.261 \\ 
SAC & 1314.8 & 1643.5 & 693.921 \\ 
FRE & 209.847 & 166.924 & 100.154 \\ 
ANA & 73.0271 & 69.5496 & 31.2973 \\ 
RIV & 182.096 & 173.425 & 78.0413 \\ 
SA & 99.6112 & 79.2361 & 47.5417 \\ 
HB & 88.0815 & 85.9332 & 40.8183 \\ 
STKN & 181.999 & 173.333 & 77.9998 \\ 
GNDL & 67.0129 & 49.0338 & 47.3993 \\ 
BKD & 499.498 & 487.315 & 231.475 \\ 
FMT & 385.036 & 366.701 & 165.015 \\ 
MOD & 315.275 & 300.262 & 135.118 \\ 
SB & 86.6804 & 82.5528 & 37.1488 \\ 
\bottomrule &  &  & 
\end{tabular}%
\end{table}

\end{document}